\documentclass{amsart}
\usepackage{mathrsfs}
\usepackage{amssymb}
\usepackage{amsfonts}
\usepackage{amsmath, color, esvect}
\usepackage[colorlinks,linkcolor=blue,urlcolor=cyan,citecolor=blue]{hyperref} 
\usepackage{genyoungtabtikz}

\usepackage[normalem]{ulem}

\usepackage{braket}
\usepackage{cleveref}

\usepackage[backend=biber,style=alphabetic,sorting=nyt]{biblatex}
\makeatletter
\def\inn#1#2{\left\langle 
      \def\ta{#1}\def\tb{#2}
      \ifx\ta\@empty{\;} \else {\ta}\fi ,
      \ifx\tb\@empty{\;} \else {\tb}\fi
      \right\rangle} 
\makeatother

\def\fhh{\mathfrak{h}}
\def\fbb{\mathfrak{b}}
\def\fnn{\mathfrak{n}}
\def\fpp{\mathfrak{p}}

\def\fbb{\mathfrak{b}}
\def\fgg{\mathfrak{g}}

\def\cC{\mathcal{C}}

\def\cO{\mathcal{O}}

\def\gr{\mathrm{gr}}
\def\bR{\mathbb{R}}
\def\fD{\mathfrak{D}}

\newcommand{\bb}{\mathbb}

\newcommand{\al}{\alpha}

\newcommand{\vf}{\varphi}

\newcommand{\lam}{\lambda}

\newcommand{\ep}{\varepsilon}

\newcommand{\fra}{\mathfrak{a}}
\newcommand{\frb}{\mathfrak{b}}

\newcommand{\frg}{\mathfrak{g}}
\newcommand{\frh}{\mathfrak{h}}

\newcommand{\frl}{\mathfrak{l}}

\newcommand{\frn}{\mathfrak{n}}

\newcommand{\frp}{\mathfrak{p}}

\newcommand{\fru}{\mathfrak{u}}

\newcommand{\caC}{\mathcal{C}}

\newcommand{\caI}{\mathcal{I}}
\newcommand{\caJ}{\mathcal{J}}

\newcommand{\caO}{\mathcal{O}}

\newcommand{\caV}{\mathcal{V}}

\newcommand{\caZ}{\mathcal{Z}}

\newcommand{\aff}{\mathbf{a}}
\newcommand{\mc}[1]{\mathcal{#1}}

\newcommand{\mf}[1]{\mathfrak{#1}}
\renewcommand{\subset}{\subseteq}
\newcommand{\rar}{\rightarrow}

\newcommand{\hs}{ \mathfrak{h}^*}

\newcommand{\Ann}{\operatorname{Ann}}
\newcommand{\gkd}{\operatorname{GKdim}}

\newtheorem{Thm}{Theorem}[section]
\newtheorem{Cor}[Thm]{Corollary}
\newtheorem{Pro}[Thm]{Proposition}
\newtheorem{Lem}[Thm]{Lemma}

\theoremstyle{definition}
\newtheorem{Rem}[Thm]{Remark}
\newtheorem{definition}[Thm]{Definition}

\newtheorem{example}[Thm]{Example}

\newtheorem{thmA}{Theorem}

\numberwithin{equation}{section}

\newcommand{\cOmin}{\cO_{\mathrm{min}}}

\newcommand{\trivial}[2][]{\if\relax\detokenize{#1}\relax{\color{red}\vspace{0em} $[$  #2 $]$}\else
\ifx#1h
\ifcsname showtrivial\endcsname{\color{orange} \vspace{0em}  $[$ #2 $]$}\fi\else{\colorWrong argument!}\fi\fi\ignorespaces}

\begin{document}

\title[Minimal highest weight modules]{Associated varieties of minimal highest weight modules}
\author{Zhanqiang Bai, Jia-Jun Ma, Wei Xiao* and Xun Xie}

\address[Bai]{School of Mathematical Sciences, Soochow University, Suzhou 215006, P. R. China}
\email{zqbai@suda.edu.cn}
\address[Ma]{Department of Mathematics, School of Mathematical Sciences, Xiamen University, Xiamen 361005, P. R. China}
\email{hoxide@gmail.com}

\address[Xiao]{School of Mathematical Sciences, Shenzhen Key Laboratory of Advanced Machine Learning and Applications, Shenzhen University,
Shenzhen 518060, Guangdong, P. R. China}
\email{xiaow@szu.edu.cn}

\address[Xie]{School of Mathematics and Statistics, Beijing
Institute of Technology, Beijing 100081, P. R. China}
\email{xieg7@163.com}



\thanks{$*$ Corresponding author.}

\long\def\mjjd#1#2{{\color{purple} #1 \sout{#2}}}


\thanks{Z. Bai was supported in part by NSFC Grant No. 12171344 and the National Key $\textrm{R}\,\&\,\textrm{D}$ Program of China (No. 2018YFA0701700 and No. 2018YFA0701701). 
J.-J. Ma was supported in part by NSFC Grant No. 11701364 and  No. 11971305, the Fundamental Research Funds for the Central Universities (Grant No. 20720230022) and  Xiamen University
Malaysia Research Fund (Grant No. XMUMRF/2022-C9/IMAT/0019).
W. Xiao was supported in part by NSFC Grant No. 11701381 and No. 12371031. X. Xie was supported in part by NSFC Grant No. 12171030 and No. 12431002.}

\begin{abstract}
	\trivial[h]{that $\mathfrak{g}$ is a complex simple Lie algebra with dual space $\frg^*$. In general, it is a difficult problem to determine the associated variety for a highest weight module of $\frg$. Let $\caO_M\subset\frg^*$ be the nilpotent coadjoint orbit associated with the annihilator of the simple $U(\frg)$-module $M$. Then $M$ is called minimal when $\caO_M$ is the unique minimal coadjoint orbit in $\frg^*$. The main result of this paper is a full classification of minimal highest weight modules for $\frg$, along with determining the associated varieties of these modules. The associated variety of a minimal highest weight module is a union of minimal orbital varieties. Since a highest weight module is considered a weak quantization of its associated variety, this amounts to giving all possible weak quantizations of minimal orbital varieties (including unions of them). Therefore, our results extend Joseph's work on minimal highest weight modules annihilated by completely prime ideals. 
 } 

	Let $\mathfrak{g}$ be a complex simple Lie algebra.  
    A simple $\fgg$-module is called minimal if the associated variety of its annihilator ideal coincides with the closure of the minimal nilpotent coadjoint orbit.    
   The main result of this paper is a classification of minimal highest weight modules for $\frg$. This classification extends the 
   work of Joseph, which focused on categorizing minimal highest weight modules annihilated by completely prime ideals.
   Furthermore, we have determined the associated varieties of these modules. In other words, we have identified all possible weak quantizations of minimal orbital varieties. 
   


\end{abstract}


\maketitle

\tableofcontents
\section{Introduction}

Let $\mathfrak{g}$ be a complex simple Lie algebra with adjoint group $G$. For every finitely generated $U(\mathfrak{g})$-module $M$, Bernstein \cite{Be} introduced the associated variety $V(M)$ in the dual space $\mathfrak{g}^{\ast}$ of $\frg$, whose dimension is equal to the Gelfand--Kirillov (GK) dimension of $M$ \cite{Vo78}. The associated variety is a useful invariant to investigate $U(\frg)$-modules, particularly highest weight modules \cite{Bai-Hu, BoB3, Ta, Mel, Mc00, Wi} and Harish-Chandra modules \cite{Vo91, Garfinkle1993, Vo19}, and has generated considerable interest in representation theory. In general, the calculation of associated varieties is difficult. Partly because the associated variety of a simple module could be reducible \cite{Ta} even in the case of type $A$ \cite{Wi}. 

If $M$ is a simple $U(\fgg)$-module, the annihilator ideal $\Ann(M)$ of $M$ is a primitive ideal of $U(\frg)$. The graded module $\mathrm{gr}(\Ann(M))$ is an ideal of the symmetric algebra $S(\frg)\simeq\mathrm{gr}(U(\frg))$.
Since $S(\frg)$ is the ring of regular functions on $\frg^*$ (see for example \cite{Har77}), the zero set $\caZ(\mathrm{gr}(\Ann(M)))$ is a subvariety of $\frg^*$, which is known to be the closure of a nilpotent coadjoint $G$-orbit $\caO_M$ of $\frg^*$ \cite{CM}. 

It is reasonable to anticipate a more elegant characterization of associated varieties for simple $U(\frg)$-modules $L$ with $\caO_L=\caO_{min}$, where $\caO_{min}$ is the unique minimal nilpotent coadjoint orbit of $\frg^*$. In this case, $L$ is referred to as \textit{minimal}. Set $\caJ=\Ann(L)$. Since $L$ is minimal, the radical $\sqrt{\mathrm{gr}\caJ}$ of the graded module $\mathrm{gr}\caJ$ is exactly the ideal $\caI(\caO_{min})$. 


The special case $\sqrt{\mathrm{gr}\caJ}=\mathrm{gr}\caJ$ (which means $\mathrm{gr}\caJ$ is prime and hence $\caJ$ is completely prime) has been extensively studied in the literature, especially when $L$ is a minimal Harish-Chandra module (see \cite{GS, BBL,  Hua99, Li00,  Sun08, HKM}). In many papers, the term ``minimal" is devoted to simple modules $L$ of this special case since $\caJ$ is now maximal among all the ideals. If $\frg$ is not of type $A$, $\caJ$ is unique and equal to the Joseph ideal $\caJ_0$ \cite{Jo98}; if $\frg$ is of type $A$, there is a family of completely prime primitive ideals $\caJ_a$ ($a\in\mathbb{C}$) \cite{Jo98, Tamori} such that $\overline{\caO}_{min}=\caZ(\mathrm{gr}\caJ_a)$.

It gets more complicated when $\sqrt{\gr\caJ} = \gr \caJ$ is not assumed. In this general situation, the minimal highest weight representations of types $A$ and $C$ are obtained by Mathieu, where they play fundamental roles in the classification of weight modules \cite{Ma}. Sun also got the minimal lowest weight modules for type $C$ when considering Harish-Chandra modules for the metaplectic double covering of $Sp_{2n}(\mathbb{R})$ \cite{Sun08}. Meng \cite{meng11} found that some unitary minimal highest weight modules satisfy a family of quadratic relations when he reconstructed the various Kepler-type problems in the unified language of Euclidean Jordan algebras. 

This paper aims to determine the associated varieties of all the minimal highest weight modules for simple Lie algebras (including the case $\caJ$ is not completely prime). Therefore, our results extend the results in \cite{BXX} for minimal highest weight Harish-Chandra modules. They can also be regarded as a first step toward classifying all the minimal representations (especially minimal Harish-Chandra modules) and determining their associated varieties.

We now state our main results. We identify $\frg^*$ with $\frg$ by the Killing form. Let $B$ be a fixed Borel subgroup of $G$ and   $\frb=\frh\oplus\frn$ be the Lie algebra of $B$ with the nilpotent radical $\frn$ and a fixed Cartan subalgebra $\frh$. 
Let $\Phi$ be the root system of $(\frg, \frh)$ with the positive system $\Phi^+$ and simple system $\Delta$ corresponding to $\frb$. Denote by $W$ the Weyl group of $\Phi$. For $\lambda\in\mathfrak{h}^*$, define 
\begin{align}\label{eq:integral_simple}
	\Phi_{[\lambda]}:&=\{\alpha\in\Phi\mid\langle\lambda, \alpha^\vee\rangle\in\mathbb{Z}\}, \\
 \label{eq:integral_Weyl}
	W_{[\lambda]}:&=\{w\in W\mid w\lambda-\lambda\in \mathbb{Z}\Phi\}.
\end{align}
Then $\Phi_{[\lambda]}$ is a root system with Weyl group $W_{[\lambda]}$ \cite[Theorem~3.4]{Hum08}. Here $\langle-,-\rangle$ denotes the natural pair between $\frh^*$ and $\frh$ and $\alpha^\vee$ represents the coroot corresponding to $\alpha$. Put 
\[
I_\lambda=\{\alpha\in\Delta\mid\langle\lambda, \alpha^\vee\rangle\in\mathbb{Z}_{>0}\}.
\]
Denote by $L(\lambda)$ the simple module with the highest weight $\lambda-\rho$, where $\rho$ is half the sum of positive roots. The associated varieties of minimal highest weight modules are given by the following result.

\smallskip

\begin{thmA}[Theorem \ref{thmav}] Let $\lambda\in \frh^*$ such that $L(\lambda)$ is minimal. Then 
\begin{enumerate}
\item the set $I_\lambda$ contains all the simple short roots\footnote{We adhere to the convention that there are no short roots in a root system of type, $A$, $D$, or $E$. };
\item the associated variety  
\[
V(L(\lambda))=\bigcup_{\alpha\in\Delta\backslash I_\lambda}\overline{Be_{\alpha}}.
\]
Here $e_\alpha$ is a nonzero root vector contained in the root space $\frg_\alpha$. 
\end{enumerate}
\end{thmA}

\smallskip
We call $\lambda$ regular if $\langle \lambda,\gamma^\vee\rangle\neq 0$ for every positive root $\gamma$. 
We call $\lambda$ semiregular if there is a unique positive root $\gamma\in\Phi^+$ satisfying $\langle\lambda, \gamma^\vee\rangle=0$. 
We call $\lambda$ integral if $\Phi_{[\lambda]} = \Phi$.  

First, we describe the classification of integral minimal highest weight modules.

Let $\ell(\cdot)$ be the \textit{length function} on $W$.
For any integral weight $\lambda\in\frh^*$, we construct the linked sequence (see Definition \ref{deflk}) $\gamma_1, \gamma_2\cdots, \gamma_m\in\Delta$ of $\lambda$ by the following process. If $|I_\lambda|=|\Delta|-1$, let $\gamma_1$ be the only simple root contained in $\Delta\backslash I_\lambda$, otherwise we terminate the process and set $m=0$. Assume now  $\gamma_{k}$ is defined. If $\Delta\backslash I_{s_{\gamma_{k}}\cdots  s_{\gamma_1}\lambda}=\{\alpha\}$ and $\ell(s_\alpha s_{\gamma_{k}}\cdots s_{\gamma_1})=k+1$, then we set $\gamma_{k+1}=\alpha$. Otherwise, we terminate the process and set $m=k$. 
In this way, we obtain the longest sequence $\gamma_1, \gamma_2\cdots, \gamma_m$ such that for any $k\leq m$, $\gamma_k$ is the unique simple root with $\langle s_{\gamma_{k-1}}\cdots s_{\gamma_1}\lambda, \gamma_k^\vee\rangle\leq0$ and $k=\ell(s_{\gamma_k}\cdots s_{\gamma_2}s_{\gamma_1})$. 
\smallskip

\begin{thmA}[Theorem \ref{thmin1}] Suppose $\lambda\in \frh^*$ is integral. Let $\gamma_1, \gamma_2\cdots, \gamma_m$ be its linked sequence. Then $L(\lambda)$ is  minimal if and only if the following conditions are satisfied:
	\begin{itemize}
		\item [(1)] $\lambda$ is regular or semiregular.
		\item [(2)] $s_{\gamma_m}\cdots s_{\gamma_2}s_{\gamma_1}\lambda$ is dominant.
		\item [(3)] $\Phi$ is of type $A$, $D$ or $E$.
	\end{itemize}
\end{thmA} 
\smallskip

Nonintegral minimal highest weight modules are classified by the following theorem.

\begin{thmA}[Theorem \ref{thmni}] Suppose $\lambda\in \frh^*$ is nonintegral. Then $L(\lambda)$ is minimal if and only if  the following conditions are satisfied:
\begin{enumerate}
	\item the weight $\lambda$ is dominant regular.
	\item the root system $\Phi$  has type $A_n$ ($n\geq1$), $B_n$($n\geq2$), $C_n$($n\geq2$), $F_4$, or $G_2$. 
	\item aligning with condition (2), the integral root system $\Phi_{[\lambda]}$ has type $A_{n-1}$, $B_1\times B_{n-1}$, $D_n$, $C_4$, or $A_2$ respectively.
\end{enumerate}
\end{thmA}

\smallskip


Following Joseph \cite{Jo98}, a highest weight module can be viewed as a weak quantization of its associated variety. Combining the above theorems, we get all possible weak quantizations of minimal orbital varieties, i.e., irreducible components of $\overline{\cOmin}\cap\fnn$. 



\medskip

The structure of this paper unfolds as follows. We begin by reviewing some preliminaries on associated varieties and $\aff$ functions in Sections 2 and 3, respectively. In section 4, we determined the associated varieties for minimal highest weight modules. In section 5 and section 6, we provide characterizations of minimal highest weight modules in the integral and nonintegral cases, respectively. 
%
%

\section{Notations and preliminaries}

%
%


Recall that $\mathfrak{g}$ is a complex simple Lie algebra with adjoint group $G$. Let $B$ be the Borel subgroup of $G$ corresponding to the Borel subalgebra $\mathfrak{b}=\frh\oplus\frn$.  Denote by $W$ the Weyl group of the root system $\Phi$ of $(\frg, \frh)$, and by $ w_0 $ the longest element of $ W $. Note that any subset $I\subset\Phi$ generates a subsystem $\Phi_I\subset\Phi$ with corresponding Weyl group $W_I\subset W$ and longest element $ w_I $. When $I\subset\Delta$, let $\mathfrak{p}_I$ be the standard parabolic subalgebra corresponding to $I$ with Levi decomposition $\mathfrak{p}_I=\mathfrak{l}_I\oplus \mathfrak{u}_I$.  We will frequently drop the subscript if there is no confusion.


\subsection{Roots and Weights}\label{subsec:RS}
For each integer $m$, the Euclidean space $\bR^m$ is equipped with the bilinear form $(\_,\_)$ such that the standard basis $\{\ep_i\}_{1\leq i\leq m}$ is orthonormal.  
We now describe the realization of an irreducible root system $\Phi:=\set{\alpha_1, \cdots, \alpha_n}$ of rank $n$ as a subset of a Euclidean space following \cite[\S12.1]{Hum78}.
When $\Phi$ has type $A_n$, $B_n$, $C_n$ or $D_n$ ($n\geq 1$), we set  
$\alpha_i:=\ep_i-\ep_{i+1}$ for $1\leq i<n$ and 
$\alpha_n := \ep_{n}-\ep_{n+1}$, $\ep_n $, $ 2\ep_n $ and $ \ep_{n-1}+\ep_n $ respectively. 
When $\Phi$ has type $E_n$ ($6 \leq n\leq 8$), we set 
$\alpha_1=\frac{1}{2}(\ep_1+\ep_8-(\ep_2+\cdots+\ep_7))$, $\alpha_2=\ep_1+\ep_2$ and $\alpha_i=\ep_{i-1}-\ep_{i-2}$ for $3\leq i\leq n$.
When $\Phi$ has type $F_4$, we set $\alpha_1=\ep_2-\ep_3$, $\alpha_2= \ep_3-\ep_4$, $\alpha_3=\ep_4$ and $\alpha_4=\frac{1}{2}(\ep_1-\ep_2-\ep_3-\ep_4)$.
When $\Phi$ has type $G_2$, we set $\alpha_1=\ep_1-\ep_2$ and $\alpha_2=-2\ep_1+\ep_2+\ep_3$.


Let $\alpha^\vee\in\frh$ be the coroot for $\alpha\in\Phi$. Then $\langle\beta, \alpha^\vee\rangle=2(\beta, \alpha)/(\alpha, \alpha)$ for $\beta\in\Phi$. The maps $\alpha^\vee\rightarrow 2\alpha/(\alpha, \alpha)$ $(\alpha\in\Phi)$ can be extended to a linear isomorphism $\vf:\frh\rightarrow\frh^*$ such that $\langle\beta, \alpha^\vee\rangle=(\beta, \vf(\alpha^\vee))$. Denote $\Phi^\vee=\{\beta^\vee\ |\ \beta\in\Phi\}$. It is called the \textit{dual system} of $\Phi$. One always has $\Phi^{\vee\vee}=\Phi$. 
If $\Phi$ is simply-laced, then $\Phi^\vee\simeq\vf(\Phi^\vee)=\Phi$. In the non simply-laced cases, we have 
$B_n^\vee\simeq \vf(B_n^\vee)=C_n$, $F_4^\vee\simeq\vf(F_4^\vee)\simeq F_4$ and $G_2^\vee\simeq\vf(G_2^\vee)\simeq G_2$.

For $\lambda\in\frh^*$, recall that $L(\lambda)$ is the simple highest weight $\mathfrak{g}$-module with highest weight $\lambda-\rho$, where $\rho$ is half the sum of positive roots. Then $L(\lambda)$ admits an central character $\chi_\lambda: Z(\fgg)\rightarrow \mathbb{C}$ such that $z\cdot v=\chi_\lambda(z)v$ for any $z\in Z(\fgg)$ and $v\in L(\lambda)$, where $Z(\fgg)$ is the center of the universal enveloping algebra $U(\fgg)$ of $\fgg$. In particular, $\chi_\lambda=\chi_\mu$ if and only if $\lambda=w\mu$ for some $w\in W$. Thus the set of central characters can be parameterized by 
\begin{equation}
\fD=\left\{\nu\in\fhh^*\left|
\begin{aligned}
&\mathrm{Re}\langle\nu,\alpha^\vee\rangle\geq0\mbox{ and whenever }\mathrm{Re}\langle\nu,\alpha^\vee\rangle=0,\\
&\mbox{ then }\mathrm{Im}\langle\nu,\alpha^\vee\rangle\geq0, \mbox{ for all }\alpha\in\Delta
\end{aligned}
\right.
\right\}.
\end{equation}
We say $L(\lambda)$ has \textit{infinitesimal character} $\nu$ when $\chi_\lambda=\chi_\nu$ for some $\nu\in \mathfrak{D}$. Put $L_w:=L(-w\rho)$ for $w\in W$. 
They form all the simple highest weight modules with infinitesimal character $\rho$.

We say $\lambda\in\frh^*$ is \emph{regular} if $\langle\lambda, \alpha^\vee\rangle\neq0$ for all $\alpha\in \Phi$, otherwise $\lambda$ is called \textit{singular}. We say $\lambda$ is \textit{dominant} (resp. \textit{anti-dominant}) if $\langle\lambda, \alpha^\vee\rangle\not\in\bb{Z}_{<0}$ (resp. $\bb{Z}_{>0}$) for all $\alpha\in \Phi$. Thus any $\nu\in\fD$ is dominant.

\subsection{GK dimensions and associated varieties}\label{subsec:GK}
In this subsection, we recall two invariants
of representations of Lie algebras after \cite{Vo78} and \cite{Vo91}. 

Let $ M $ be a finite generated $U(\frg)$-module. Choose a finite dimensional generating subspace $M_0$ of $M$. Let $\{U_{n}(\mathfrak{g})\}_{n\geq 0}$ be the standard filtration of $U(\mathfrak{g})$. Write $M_n=U_n(\mathfrak{g})\cdot M_0$ and denote $M_{-1}=0$. Set
\[
\mathrm{gr} (M)=\bigoplus\limits_{n=0}^{\infty} \mathrm{gr}_n M,
\]
where $\mathrm{gr}_n M=M_n/{M_{n-1}}$. Then $\mathrm{gr}(M)$ is a graded module of $\mathrm{gr}(U(\mathfrak{g}))\simeq S(\frg)$. 


\begin{definition} The \textit{Gelfand--Kirillov dimension} of $M$  is defined by$$
	\operatorname{GKdim} M = \overline{\lim\limits_{n\rightarrow \infty}}\frac{\log\dim( U_n(\mathfrak{g})M_{0} )}{\log n}.
	$$
	The  \textit{associated variety} of $M$ is defined by
	\[
	V(M):=\{X\in \mathfrak{g}^* \mid p(X)=0 \text{ for all~} p\in \operatorname{Ann}_{S(\mathfrak{g})}(\mathrm{gr} M)\}.
	\]
\end{definition}

In particular, we have (see e.g., \cite{NOT}) 
\begin{equation}\label{eq:gkd-v}
 \dim V(M)=\gkd M.
\end{equation} 
Since we can identify $\frg^*$ with $\frg$ via the Killing form on $\frg$, $V(M)$ can also be viewed as a subvariety of $\frg$. 

One can also define the associated variety for a two-sided ideal of $U(\frg)$.

\begin{definition}
	Let $\caI$ be a two-sided ideal in $U(\mathfrak{g})$. Then $\mathrm{gr}(U(\mathfrak{g})/\caI)\simeq S(\mathfrak{g})/\mathrm{gr}\caI$ is a graded $S(\mathfrak{g})$-module with annihilator $\mathrm{gr}\caI$. The \textit{associated variety} of $\caI$ is
	\[
	V(\caI):=V(U(\mathfrak{g})/\caI)=\{X\in \mathfrak{g}^* \mid p(X)=0\ \mbox{for all $p \in {\mathrm{gr}}\caI$}\}.
	\]
\end{definition}
For a highest weight module $M$, we have (see \cite[Cor. 2.8]{Jo78})
	\begin{equation}\label{eq:dim-ass}
	\dim V(\operatorname{Ann}(M))=2\dim V(M).
	\end{equation}

If $M_0$ is $\fra$-invariant for a subalgebra $\fra$ of $\frg$, then (see \cite[(1.5)(b)]{Vo78})
\begin{equation}\label{embed}
V(M)\subset V(\operatorname{Ann}(M))\cap (\frg/\fra)^*.
\end{equation}
In particular, if $M$ is generated by a finite dimensional $\frb$-invariant space, then $V(M)\subset(\bar\frn)^*\simeq \frn$, where the last isomorphism is induced from the Killing form. 

\begin{Lem}\label{lem:vanno}
For every $\lambda\in\frh^*$, the associated variety $V(\Ann(L(\lambda)))$ is the closure $\overline{\mathcal{O}}$ of a nilpotent orbit $\mathcal{O}$.
\end{Lem}

This result can be found in \cite[3.10]{Jo85} or \cite[Theorem 10.2.2]{CM}. Suppose that $\mathcal{O}\subseteq \mathfrak{g}\simeq\frg^*$ is a nilpotent $G$-orbit with closure $\overline{\mathcal{O}}$ \cite{CM}. Each irreducible component of $\overline{\mathcal{O}}\cap \mathfrak{n}$ is called an {\it orbital variety}, which can be written as $\mathcal{V}(w):=\overline{B(\mathfrak{n}\cap w\mathfrak{n})}$ for some $w\in W$ \cite{Ta}. 

The following result is known as Spaltenstein-Steinberg equality \cite[3.1]{Jo84}.

\begin{Lem}\label{dim}
	Let $\mathbb{O}$ be an orbital variety of a nilpotent orbit $\mathcal{O}$. Then $\dim \mathbb{O}=\frac{1}{2}\dim\mathcal{O}$.
\end{Lem}

The next result is due to Joseph \cite[Proposition 4.7]{Jo84} 
and Borho-Brylinski \cite[Proposition~6.11 (a)]{BoB3}.

\begin{Lem}\label{ob}
	Let $\lambda\in\frh^*$. Then $V(L(\lambda))$ is a union of some orbital varieties contained in $\overline{\mathcal{O}}\cap\frn$, where $\overline{\mathcal{O}}=V(\Ann(L(\lambda)))$. Moreover, $\caV(w)\subset V(L_w)$.
\end{Lem}

\subsection{Minimal orbits and minimal modules}

For any simple Lie algebra $\mathfrak{g}$, there is a unique nonzero nilpotent orbit $\mathcal{O}_{min}$ of minimal dimension \cite{CM}. Then each irreducible component of  $\overline{\caO}_{min}\cap \mf{n}$ is called a \textit{minimal orbital variety} by Joseph \cite{Jo98}. 

\begin{definition}\label{defmin1}
We call a simple $U(\frg)$-module $M$ \textit{minimal} if 
\begin{equation}\label{dm1eq1}
V(\operatorname{Ann} (M))=\overline{\mathcal{O}}_{\min}.
\end{equation}
\end{definition}

\begin{Rem}
 In \cite{GS}, a Harish-Chandra module $M$ is called minimal if (\ref{dm1eq1}) holds and $\operatorname{Ann}(M)$ is a Joseph ideal, which is completely prime. In our Definition \ref{defmin1}, the annihilator ideal $\operatorname{Ann}(M)$ is not necessarily completely prime. 
\end{Rem}

The following result is an easy consequence of Lemma \ref{ob}.

\begin{Lem}\label{min}
	If $L(\lambda)$ is a minimal highest weight module, then $V(L(\lambda))$ is a union of some minimal orbital varieties.
\end{Lem}

An orbital variety is called {\it weakly quantizable} \cite{Jo98} if it is the associated variety of a highest weight module. In this paper, we also say a union of orbital varieties is {\it weakly quantizable} when it is the associated variety of a highest weight module. By Joseph \cite[\S4.12]{Jo98},  the following proposition holds.

\begin{Pro}\label{J}
	Let  $\mathbb{O}$ be a minimal orbital variety. Then
	\begin{itemize}
		\item [(1)] $\mathbb{O}=\overline{Be_{\alpha}}=\mathcal{V}(s_{\alpha}w_0)$ for some long simple root $\alpha\in\Delta$.
		\item [(2)] $\mathbb{O}$ is weakly quantizable except when $\mathbb{O}=\overline{Be_{\alpha_i}}$ for $2\leq i\leq n-2$ in type $B_n (n\geq 4)$.
	\end{itemize}
\end{Pro}

We will see that minimal modules always exist for any fixed simple Lie algebra $\fgg$. Hence, according to \eqref{eq:gkd-v}, \eqref{eq:dim-ass} and Lemma \ref{lem:vanno}, the minimal highest weight modules can be characterized as those $ L(\lambda) $ with the minimal (nonzero) GK dimension, which is equal to $\frac12\dim {\mathcal {O}}_{\min}=\langle\rho, \beta^{\vee}\rangle $ by \cite{Wang}, where $\beta\in\Phi^+$ is the highest root (e.g., \cite[\S 13.4]{Hum78}). The strategy for finding them is clear: we investigate the highest weight modules with minimal possible GK dimension when an infinitesimal character is fixed. The following notation is useful.

\begin{definition}\label{defnu}
	Fix $\nu\in\mathfrak{D}$. We say $ L(\lambda) $ is \textit{$\nu$-minimal} if $\chi_\lambda=\chi_\nu$ and $L(\lambda)$ has minimal possible nonzero GK dimension among all the simple highest weight modules with infinitesimal character $\nu$.
\end{definition}

Obviously, a minimal highest weight module is always $\nu$-minimal for some $\nu\in\mathfrak{D}$, but not vice versa. The values of the minimal and $\rho$-minimal GK dimension for each type are given in Table \ref{tab1}, see \cite{Vo81, Lus84, BXX}.

\begin{table}[htbp]
	\centering
	\renewcommand{\arraystretch}{1.5}
	\setlength\tabcolsep{5pt}
	
	\begin{tabular}{c|ccccccccc}
		\hline & $A_{n}$ & $B_n$ & $C_n$ & $ D_n $ & $ E_6 $ & $ E_7 $ & $ E_8 $ & $ F_4 $ & $ G_2 $ \\ \hline
		$ \min \text{ of }\gkd L(\lambda) $    &  $ n $  &  $ 2n-2 $     &   $ n $    &    $ 2n-3 $     &      11   &    17     &    29     &     8    &  3       \\ \hline
		$ \min \text{ of }\gkd L_w $    &  $ n $  &    $ 2n-1 $  &   $ 2n-1 $   &    $ 2n-3 $     &      11   &    17     &    29     &     11    &  5       \\ \hline
	\end{tabular}
	\bigskip
	\caption{Minimal (nonzero) GK dimensions}
	\label{tab1}
\end{table}


Minimal highest weight modules could be integral or nonintegral. For the integral case, we need another small orbit $\mathcal{O}_s$, which is called \textit{minimal special} since it is minimal among all non-trivial special nilpotent orbits (see \cite[\S6.3]{CM} for the definition of special orbit). From Humphreys \cite[Proposition~5]{Hum16}, we have $\dim \mathcal{O}_{s}=2\langle\rho, \beta_s^{\vee}\rangle $, where $\beta_s$ is the highest short root  when $\Phi$ is not simply-laced and is the highest root $\beta$ otherwise. Thus $\mathcal{O}_s=\mathcal{O}_{min}$ if and only if $\Phi$ is simply-laced. 

\begin{definition}\label{defsmin}
We call a simple $U(\frg)$-module $M$  {\it  integral minimal} if 
it has integral infinitesimal character and 
\begin{equation}
V(\operatorname{Ann} (M))=\overline{\mathcal{O}}_{s}.
\end{equation}
\end{definition}
Since  $\mathcal{O}_{s}$ is the unique nilpotent orbit with dimension $2\langle\rho, \beta_s^{\vee}\rangle$ (see \cite[Remark 3.(b)]{Lu.Green}),  a simple highest weight module $M$ with integral infinitesimal character is integral minimal if and only if 
\begin{equation}
\gkd M=\frac12\dim {\mathcal {O}}_{s}=\langle\rho, \beta_s^{\vee}\rangle 
\end{equation}
in view of (\ref{eq:dim-ass}). 
We will present a characterization of integral minimal special highest weight modules in \S 5.
%
%

\section{Algorithm for GK dimensions}

%
%

In this section, we will describe the algorithm obtained in \cite{BX} for GK dimension of highest weight modules. 

Recall $W_{[\lambda]}$ and $\Phi_{[\lambda]}$ defined in \eqref{eq:integral_Weyl} and \eqref{eq:integral_simple} respectively. Let $\Delta_{[\lambda]}$ be the simple system of $\Phi_{[\lambda]}$ contained in the set $\Phi^+$ of positive roots. For any $\lambda\in\frh^*$, there exists a unique shortest element $w_\lambda\in W_{[\lambda]}$ such that $\mu=w_\lambda^{-1}\lambda$ is anti-dominant. 
Thus $W_{[\lambda]}=W_{[\mu]}$ and $w_\lambda$ is the minimal length representative of a coset of $W_{[\mu]}/W_{J_\mu}$, where $W_{J_\mu}$ is the Weyl group of the subsystem $\Phi_{J_\mu}$ generated by 
\begin{equation}\label{eq:Jmu}
J_\mu=\{\alpha\in\Delta_{[\mu]}\mid\langle\mu, \alpha^\vee\rangle=0\}. 
\end{equation}
More results about integral subsystems can be found in \cite[\S3.4, \S3.5]{Hum08}, \cite[\S5.3]{BXX} and \cite[\S4.1]{XW}.


\begin{Pro}[{\cite[Prop. 3.8]{BX}}]\label{pr:main1}
	Let $ \lam\in\hs $. Then
	\begin{equation*}
		\gkd L(\lambda)=|\Phi^+|-\aff_{[\lambda]}(w_\lambda),
	\end{equation*}
	where $\aff_{[\lambda]}: W_{[\lambda]}\rightarrow\mathbb{N}$ is Lusztig's $\aff$-function on the integral Weyl group $W_{[\lambda]}$. 
\end{Pro}

 We now explain more about Lusztig's $\aff$-function,  $\aff: W\rightarrow \mathbb{N}$ in the following. A good reference is \cite{Lus03}. 


Recall that the Weyl group $ W  $ is a Coxeter group generated by $ S=\{s_\al\mid\al\in\Delta \} $. Given an indeterminate $v$, the Hecke algebra $ \mc{H} $ over $ \mathcal{A} :=\mathbb{Z}[q,q^{-1}]$ is generated by $ T_w $, $ w\in W $ with relations \[
T_wT_{w'}=T_{ww'} \text{ if }\ell(ww')=\ell(w)+\ell(w'),
\]
\[
\text{and }(T_s+q^{-1})(T_s-q)=0 \text{ for any }s\in S.
\]
The unique elements $ C_w $ such that
\[
\overline{C_w}=C_w,\qquad C_w\equiv T_w \mod{\mc{H}_{<0}}
\]
are known as the \textit{Kazhdan-Lusztig} (KL) \textit{basis} of $ \mc{H} $, where $ \bar{\,} :\mc{H}\rar\mc{H}$ is the bar involution such that $ \bar{q}=q^{-1} $, $ \overline{T_w} =T_{w^{-1}}^{-1}$, and $ \mc{H}_{<0}=\bigoplus_{w\in W}\mathcal{A}_{<0}T_w $ with $ \mathcal{A}_{<0}=q^{-1}\mathbb{Z}[q^{-1}] $.

If $ C_y $ occurs in the expansion of $ hC_w $ (resp. $C_wh$) with respect to the KL basis for some $ h\in\mc{H} $, then we write $ y\leftarrow_L w $ (resp. $ y\leftarrow_R w $). Extend $ \leftarrow_L $ (resp. $ \leftarrow_R $) to a preorder $ \leq_L $ (resp. $\leq _R$) on $ W $. For $x, w\in W$, write $x \leq_{LR} w$ if there exists $x=w_1, \cdots, w_n=w$ such that for every $1\leq i<n$ we have either $w_i\leq_L w_{i+1}$ or $w_i\leq_R w_{i+1}$. Let $\sim_{L}$, $\sim_{R}$, $\sim_{LR}$ be the equivalence relations associated with $\leq_L$, $\leq_R$, $\leq_{LR}$ (for example, $x\sim_{L}w$ if and only if $x\leq_L w$ and $w\leq_Lx$). The equivalence classes on $W$ for $\sim_L$, $\sim_R$, $\sim_{LR}$ are called \textit{left cells}, \textit{right cells} and \textit{two-sided cells} respectively. The following result is straightforward consequence of \cite[(2.3a) and (2.3b)]{KL}.

\begin{Lem}\label{hklem1}
	Let $x, y, z\in W$. If $z=xy$ and $\ell(z)=\ell(x)+\ell(y)$, then $z\leq_L y$ and $z\leq_R x$.
\end{Lem}

Obviously $\{1\}$ (resp. $\{w_0\}$) is the largest (smallest) two sided-cell of $W$. Let 
\[
\cC := \Set{w\in W | \text{$w\neq 1$ and has a unique   reduced expression}}.
\]
Evidently, $S\subset\caC$.
Denote \[
\mathcal{C}w_0:=\{ww_0\mid w\in\mathcal{C}\}. 
\]

\begin{Lem}\label{hklem2}
Let $w\in W$. Then
	\begin{itemize}
 \item[(1)]  The sets $\mathcal{C}$ and $\mathcal{C}w_0$ are two-sided cells of $W$.
		\item[(2)]  If $w\neq 1, w_0$, then $x\leq_{LR} w\leq_{LR} y$ for any $x\in\mathcal{C}w_0$ and $y\in\mathcal{C}$.
            \item[(3)] The equations $\cC = w_0^{-1} \cC w_0$ and $\cC w_0 = w_0\cC$ hold.
	\end{itemize}

\end{Lem}
\begin{proof}
Part (1) is \cite[Prop. 3.8~(c)]{Lus83} combining with \cite[Remark 3.3~a)]{KL}. For part (2), since $w\neq 1$, there exists $s\in S$ such that $\ell(sw)<\ell(w)$. Lemma \ref{hklem1} implies that $w\leq_R s$. Thus $w\leq_{LR} y$ in view of $s\in S\subset\mathcal{C}$ and (1). Similarly, we can show that $x\leq_{LR} w$.
Part (3) follows from the fact that $s \mapsto w_0^{-1} s w_0$ induces an involution on the set of simple roots. 
\end{proof}

Write $ C_xC_y=\sum_{z\in W} h_{x,y,z}C_z $ with $ h_{x,y,x}\in\mathcal{A} $. Then 
$ \mathbf{a}:W\rar\mathbb{N} $ is defined by\[
\aff(z)=\max\{\deg h_{x,y,z}\mid x,y\in W \} \text{ for } z\in W.
\]

The following lemma is an easy consequence of Lusztig's results \cite[\S15]{Lus03}.

\begin{Lem}\label{alem1} 
	Let $x, w\in W$. Then
	\begin{itemize}\item [(1)]$ \aff(w)=\aff(w^{-1}) $.
		\item [(2)] If $x\leq_{LR} w$, then $\aff(x)\geq\aff(w)$. Hence $\aff(x)=\aff(w)$ whenever $x\sim_{LR} w$.
		\item [(3)] If $ w_I $ is the longest element of the parabolic subgroup of $ W $ generated by simply reflections of $ I\subset \Delta $, then $ \aff(w_I)$ is equal to  the length $\ell(w_I) $ of $ w_I $.
		\item [(4)] If $ W $ is a direct product of  Coxeter subgroups $ W_1 $ and $ W_2 $, then\[
		\aff(w)=\aff(w_1)+\aff(w_2)
		\]for  $ w=(w_1,w_2) \in W_1\times W_2=W$.
	\end{itemize}
\end{Lem}

\begin{Lem}\label{alem2}
	We have the following equation 
	\[
	\aff(1)=0<1=\aff(\mathcal{C})<\aff(w)<\aff(\mathcal{C}w_0)<\aff(w_0)=\ell(w_0)=|\Phi^+|.
	\]
       with $w\in  W \setminus  (\mathcal{C}\cup\mathcal{C}w_0\cup\{1, w_0\})$.
\end{Lem}
\begin{proof}
	In view of \cite[Prop. 13.7, 13.8]{Lus03}, we get $\aff(1)=0<\aff(\mathcal{C})$ and $\aff(\mathcal{C}w_0)<\ell(w_0)=\aff(w_0)=\ell(w_0)$. Lemma \ref{alem1} implies $\aff(\mathcal{C})=\aff(s)=1$ for any $s\in S$, while Lemma \ref{hklem2} yields $\aff(\mathcal{C})\leq\aff(w)\leq\aff(\mathcal{C}w_0)$. If $\aff(\mathcal{C})=\aff(w)$, we obtain $w\in\mathcal{C}$ by \cite[\S14.2 P11]{Lus03} and \cite[Corollary 2.15]{BV83}, a contradiction. We can get a similar contradiction when we assume $\aff(w)=\aff(\mathcal{C}w_0)$.
\end{proof}

\begin{Rem}\label{smRm}
    From Proposition \ref{pr:main1} and Lemma \ref{alem2}, an integral highest weight module $M$ is integral minimal if and only if 
\begin{equation}\label{eqmi1}
\gkd M=\frac12\dim {\mathcal {O}}_{s}=\langle\rho, \beta_s^{\vee}\rangle
=|\Phi^+|-\aff(w_0\mathcal{C}). 
\end{equation}
Thus $\aff(w_0\mathcal{C})
=|\Phi^+|-\langle\rho, \beta_s^{\vee}\rangle$.
\end{Rem}
%
%

\section{Associated varieties of highest weight modules}

%
%
In this section, we determine the associated varieties of all the minimal highest weight modules (Theorem \ref{thmav}). We will also recall the relation between associated varieties and cells.

\subsection{Associated varieties of minimal highest weight modules} Note that any subset $I\subseteq\Delta$ generates a subsystem $\Phi_I\subseteq\Phi$ with corresponding parabolic subalgebra $\mathfrak{p}=\mathfrak{l}\oplus \mathfrak{u}\supseteq\frb$ of $\frg$.

We say $\lambda\in\frh^*$ is \emph{$\Phi_I^+$-dominant} if and only if $\langle\lambda, \alpha^\vee\rangle\in\mathbb{Z}_{>0}$ for all $\alpha\in I$. Let $\bar\fru$ be the nilpotent subalgebra opposite to $\fru$. Then $\bar\frn = \bar\fru \oplus \bar\frn_I$ and $\fpp = \bar\frn_I \oplus \fbb$ where $\bar\frn_I := \bar\frn \cap \frl$. The following result is a generalization of \cite[Lem. 6.1]{Jo84}.

\begin{Pro}\label{HC}
	Let $\frg$ be a simple complex Lie algebra and $\lambda\in \frh^*$. Then $\lambda$ is $\Phi_I^+$-dominant if and only if $V(L(\lambda))\subset\fru$.
\end{Pro}

\begin{proof}
	Denote by $F(\lambda)$ the simple $\frl$-module with highest weight $\lambda-\rho$. Then $F(\lambda)$ is finite dimensional if and only if $\lambda$ is $\Phi_I^+$-dominant.
	
    Assume that $\lambda$ is $\Phi_I^+$-dominant. Since $F(\lambda)$ can be viewed as a $\frp$-module with trivial $\fru$-action, we obtain a surjective homomorphism
	\[
	\varphi: U(\frg)\otimes_{U(\frp)}F(\lambda)\rightarrow L(\lambda).
	\]
	Then $L(\lambda)$ is generated by $\vf(F(\lambda))\simeq F(\lambda)$, which is $\frp$-invariant. Thus (\ref{embed}) implies that $V(L(\lambda))\subset(\frg/\frp)^*\simeq(\bar\fru)^*\simeq\fru$.
	
 Now assume that $V(L(\lambda))\subset\fru$.
 Let $v^+$ be a nonzero highest weight vector of $M:=L(\lambda)$.
 Let $M_0:=\mathbb{C}v^+$, which is an one dimensional $\frb$-invariant generating space of $M$.
 Consider the $\frb$-invariant filtration $M_n:=U_n(\frg)\cdot M_0=U_n(\bar\frn)\cdot M_0$.
 Since the subvariety $(\bar\fru)^* = \fru \subset\frn$ is defined by the ideal $\bar\frn_IS(\bar\frn)$, the condition $V(M)\subset \fru$ means that the $\bar\frn_I$-action on $\gr M$ is nilpotent.
 In other words, there is a positive integer $j$ such that $U_j(\bar\frn_I)\cdot M_i \subset M_{i+j-1}$ for all $i\in \mathbb{N}$.
 For every $\alpha\in I$, let $e_{-\alpha}$ be a non-zero vector in the $-\alpha$-root space $\fgg_{-\alpha}$.  Then we have   
    \begin{equation}\label{HCeq1}
		e_{-\alpha}^j \cdot v^+=u\cdot v^+
    \end{equation}
    for some $u\in U_{j-1}(\bar{\frn})$. 
    We claim that $u\cdot v^+ = 0$. Otherwise, we can assume that $u$ has weight $j\alpha$ by comparing the weights on the two sides of \eqref{HCeq1}.  
    Since $u\in U_{j-1}(\bar\fru)$, there must be a sequence $\beta_1, \beta_2, \cdots, \beta_k$  in $\Phi^+$ for some  $k\leq j-1$ such that    
     \[
     -(\beta_1+\beta_2+\cdots+\beta_k) = -j\alpha . 
     \]
     But this can not hold since $\alpha$ is a simple root.
     Hence $e_{-\alpha}^j\cdot v^+=0$. 
     By $\mathfrak{sl}_{2}$ theory, we conclude that $\langle\lambda, \alpha^{\vee}\rangle\in \mathbb{Z}_{>0}$ for every $\alpha\in I$, i.e., $\lambda$ is $\Phi_I^+$ dominant.
\end{proof}


Recall that $I_\lambda=\{\alpha\in\Delta\mid\langle\lambda, \alpha^\vee\rangle\in\mathbb{Z}_{>0}\}$. Obviously $\lambda$ is $\Phi_{I_\lambda}^+$-dominant.

\begin{Thm}\label{thmav}
	Let $\lambda\in \frh^*$. If $L(\lambda)$ is a minimal highest weight module, then
	\[
	V(L(\lambda))=\bigcup_{\alpha\in\Delta\backslash I_\lambda}\overline{Be_{\alpha}}.
	\]
	Moreover, $I_\lambda$ contains all the simple short roots.
\end{Thm}
\begin{proof}
	Set $I=I_\lambda$. By Lemma \ref{min} and  Proposition \ref{J}, we can write $V(L(\lambda))=\cup_{\alpha\in T}\overline{Be_{\alpha}}$ for a subset $T\subset\Delta$ which contains no short root. Since $\lambda$ is $\Phi_{I}^+$-dominant, Proposition \ref{HC} yields $e_\alpha\in\fru_I$ for $\alpha\in T$. This forces $T\subset\Delta\backslash I$. On the other hand, with $e_{\alpha}\in\fru_{\Delta\backslash T} $ for $\alpha\in T$, one has $V(L(\lambda))\subset\fru_{\Delta\backslash T}$. Thus $\lambda$ is $\Phi_{\Delta\backslash T}^+$-dominant by Proposition \ref{HC}. We obtain $\Delta\backslash T\subset I$. Hence $T=\Delta\backslash I$.
\end{proof}

\trivial[h]{

\subsection{Right cells and associated varieties}

Denote $\caI_w=\Ann(L_w)$. By Borho-Brylinski \cite{BoB1} and Joseph \cite{Jo85}, $V(\caI_w)=V(U(\mathfrak{g})/\caI_w)$ is the closure of a single special (in the sense of Lusztig \cite{Lu79}) nilpotent orbit. From \cite{KL} and \cite{BV83}, there is a bijection between special nilpotent orbits of $\frg$ and two-sided cells of the Weyl group $W$. \mjjd{Delete}{, see also \cite{Ta}.
Some more properties of special orbits can be found in Collingwood-McGovern \cite{CM}. }

\mjjd{Move above}{For a nilpotent orbit $\mathcal{O}$, denote by $\mathrm{Irr}(\overline{\mathcal{O}}\cap \mathfrak{n})$ the set of irreducible components in $\overline{\mathcal{O}}\cap \mathfrak{n}$.}
The following result was conjectured by Tanisaki \cite[Conj.3.4]{Ta}. The case of type $A$ was proved by Borho-Brylinski \cite{BoB3}, while the case of the \mjjd{type $B$, $C$ or $D$}{other classical types} was proved by McGovern \cite[Section~3]{Mc00}.

\begin{Thm}\label{Ta}
	Suppose that $\frg$ is of type $A$, $B$, $C$ or $D$. \mjjd{}{classical.} Let  $\mathscr{C}$ be the two-sided cell corresponding to a special nilpotent orbit $\mathcal{O}$. Denote by $\mathscr{C}_R$  the set of right cells contained in $\mathscr{C}$. Then there exists a bijection from $\mathscr{C}_R$ to $\mathrm{Irr}(\overline{\mathcal{O}}\cap \mathfrak{n})$ $(w\rightarrow Y_w)$ and an ordering $\prec$ on $\mathrm{Irr}(\overline{\mathcal{O}}\cap \mathfrak{n})$ such that $V(L_w)=Y_w \cup \tilde Y_w$, where $\tilde{Y}_w$ is a union of some orbital varieties $Y$ in $\mathrm{Irr}(\overline{\mathcal{O}}\cap \mathfrak{n})$  with $Y\prec Y_w$.
\end{Thm}

\begin{Rem}
	In Tanisaki's original conjecture \cite[Conj.3.4]{Ta}, $\frg$ is not necessarily classical, where the cases of type $E_6$ and $G_2$ were proved by himself. For the case of type $F_4$, Tanisaki showed that the conjecture is true for nine special nilpotent orbits among all the eleven ones. The other cases are still unknown.
\end{Rem}

\begin{Pro}\label{m1}
	Suppose $\mathfrak{g}$ is of \mjjd{}{classical type}. Let  $w,y\in W$, then $w\sim_R y$ if and only if $V(L_w)=V(L_y)$.
\end{Pro}

\begin{proof}
	If $w\sim_R y$, we get $V(L_w)=V(L_y)$ by \cite[Lem. 6.6]{Jo84} and \cite[Cor. 6.3]{BoB3}. Now assume that $V(L_w)=V(L_y)$. In view of \cite[Cor. 4.11]{BoB1} and \cite[Thm. 10.2.2]{CM}, $V(\caI_w)=V(\caI_y)$ is the closure of a special nilpotent orbit $\mathcal{O}$. Let $\mathscr{C} $ be the two-sided cell corresponding to $\mathcal{O}$ in the Springer correspondence. With $V(L_w)=V(L_y)$, Theorem \ref{Ta} yields $Y_w \cup \tilde{Y}_w=Y_y \cup \tilde{Y}_y$. This forces $Y_w=Y_y$. Hence $w\sim_R y$. 
	
\end{proof}

\begin{Rem}\label{EG} The proposition holds for all the simple Lie algebras once Tanisaki's original conjecture \cite[Conj. 3.4]{Ta} is proved completely.
\end{Rem}

From \cite[Lem. 5.2]{BBM}, we know that the associated variety of a simple $\mathfrak{g}$-module is invariant under the corresponding translation functor. Recall that $\mathcal{V}(w):=\overline{B(\mathfrak{n}\cap w\mathfrak{n})}$ is the orbital variety for $w\in W$.  For any $\lambda\in\frh^*$, recall that there exists a unique shortest element $w_\lambda\in W_{[\lambda]}$ such that $\mu=w_\lambda^{-1}\lambda$ is anti-dominant. 
Combined with Lemma \ref{ob} and \cite[Cor. 3.3]{BX}, we have the following results.

\begin{Cor}\label{l3.3}
	Let $\lambda\in\frh^*$ be an integral weight. Then $$ V(L(\lambda))=V(L_{w_\lambda})\supseteq\caV(w_\lambda) .$$
\end{Cor}


\begin{Cor}\label{thm:rightcell}
	Suppose that $\frg$ is classical. Let $\lambda, \mu\in\frh^*$ be integral. Then $V(L(\lambda))=V(L(\mu))$ if and only if $w_\lambda\sim_Rw_\mu$.
\end{Cor}

}
%
%

\section{Integral minimal highest weight modules}

%
%
Recall that an integral minimal highest weight module must be minimal special. In this section, we will provide the description of integral minimal highest weight modules. The following definition is useful.

\begin{definition}\label{deflk}
	Let $\lambda\in\frh^*$ be integral. Choose the longest sequence $\gamma_1, \gamma_2\cdots, \gamma_m\in\Delta$ such that for any $k\leq m$, $\gamma_k$ is the unique simple root with $\langle s_{\gamma_{k-1}}\cdots s_{\gamma_1}\lambda, \gamma_k^\vee\rangle\leq0$ $(s_{\gamma_{0}}=1)$ and $k=\ell(s_{\gamma_k}\cdots s_{\gamma_2}s_{\gamma_1})$. We call $\gamma_1, \gamma_2\cdots, \gamma_m$ the \textit{linked sequence} of $\lambda$.
\end{definition}

By the construction of linked sequence,  
$s_{\gamma_k}\cdots s_{\gamma_2}s_{\gamma_1}$ is a reduced expression and in $\cC$ for $1 \leq k\leq m$.

\begin{Thm}\label{thmin1}
	Let $\lambda\in\frh^*$ be integral and $\gamma_1, \gamma_2\cdots, \gamma_m$ be the linked sequence of $\lambda$. Then the following conditions are equivalent:
	\begin{enumerate}
		\item The highest weight module $L(\lambda)$ is integral minimal.
		\item The number $m$ is positive and  $w_0w_\lambda^{-1}=s_{\gamma_m}\cdots s_{\gamma_2}s_{\gamma_1}$.  
		\item  The number $m$ is positive, $s_{\gamma_m}\cdots s_{\gamma_2}s_{\gamma_1}\lambda$ is dominant, and $\lambda$ is regular or semiregular.
	\end{enumerate}
\end{Thm}
Recall that we say an integral weight $\lambda\in\frh^*$ is \textit{semiregular} if there is a unique positive root $\gamma\in\Phi^+$ such that $\langle\lambda, \gamma^\vee\rangle=0$.

\begin{proof}
	First, we show that (1) implies (2).
 If $L(\lambda)$ is minimal special, then by Remark \ref{smRm} we have \[\gkd L(\lambda)=|\Phi^+|-\aff(w_0\mathcal{C}).\]
 Set $w=w_0w_\lambda^{-1}$.
 It is the longest element such that $\nu=w\lambda$ is dominant.
 Proposition \ref{pr:main1} and Lemma \ref{alem2} yield $w_\lambda^{-1}\in w_0\mathcal{C}$.
 Therefore $w=w_0w_\lambda^{-1}\in\mathcal{C}$ has a unique reduced expression.
 So we can find a unique $\alpha_1\in\Delta$ such that $\ell(ws_{\alpha_1})<\ell(w)$, i.e., $w\alpha<0$.
 One has $0\geq\langle w\lambda, (w\alpha_1)^\vee\rangle=\langle \lambda, \alpha_1^\vee\rangle$.  Thus $1\leq m$  $\gamma_1 = \alpha_1$ under the notation of \Cref{deflk}.
 Similarly, if $ws_{\alpha_1}\neq 1$, then there is a unique $\alpha_2 \in \Delta$ such that $\ell(ws_{\alpha_1}s_{\alpha_2})<\ell(ws_{\alpha_1})$. Moreover, $2\leq m$ and  $\alpha_2 = \gamma_2$. Inductively,  can eventually get certain $t \leq m$ and  $\alpha_1, \cdots, \alpha_t \in \Delta$ such that  $ws_{\alpha_1}s_{\alpha_2}\cdots s_{\alpha_t}=1$ and $\alpha_i = \gamma_i$ for $1\leq i \leq t$.
 We claim that $t=m$. Otherwise, $t<m$ and then  
	$\langle w\lambda, \gamma^\vee_{t+1}\rangle=\langle s_{\gamma_t}\cdots s_{\gamma_1}\lambda, \gamma^\vee_{t+1}\rangle\leq0$. Since $w\lambda$ is dominant, this forces 
	$\langle w\lambda, \gamma^\vee_{t+1}\rangle=0$. Therefore $(s_{\gamma_{t+1}}w)\lambda=w\lambda$ is dominant and $\ell(s_{\gamma_{t+1}}w)>\ell(w)$, this contradict to the fact that $w$ is the longest element such that $w\lambda$ is dominant. 
 Hence we conclude that $s_{\gamma_m}\cdots s_{\gamma_1} = w\in\caC$ and is a reduced expression. 
	
  Next, we show that (2) implies (3).
Since $w=w_0w_\lambda^{-1}=s_{\gamma_m}\cdots s_{\gamma_1}\in \caC$ is the longest element such that $w\lambda$ is dominant, it suffices to prove $\lambda$ is either regular or semiregular.
Denote $I=\{\alpha\in\Delta\mid \langle w\lambda, \alpha^\vee\rangle=0\}$.
For every $\alpha\in I$,  one must have $\ell(s_\alpha w)<\ell(w)$ (so  $s_\alpha$ occurs on the left of some reduced expression of $w$)  since $s_\alpha w\lambda=w\lambda$ is also dominant.  
Since $w$ has a unique reduced expression, we conclude that $|I|\leq 1$, which implies $\lambda$ is regular or semiregular. 
	
  Now we show that (3) implies (1).
  First, $x=s_{\gamma_m}\cdots s_{\gamma_1}$ has the unique reduced expression by the construction in \Cref{deflk}, that is, $x\in\caC$.
  Moreover, $w\lambda=x\lambda$ since both sides are dominant weights in $W\lambda$.
  If $\lambda$ is regular, we have $w=x$.
  We obtain $w_\lambda=x^{-1}w_0\in\caC w_0$ and $L(\lambda)$ is integral minimal, in view of Proposition \ref{pr:main1} and Lemma \ref{alem2}.
  Now suppose that $\lambda$ is semiregular. Then $\langle x\lambda, \alpha^\vee\rangle=0$ for a unique $\alpha\in\Delta$. Note that $w\lambda=x\lambda=s_\alpha x\lambda$.
  If $\ell (s_\alpha x)<\ell(x)$, then $w=x$ and the remaining argument is the same as the regular case.
  The case $\ell (s_\alpha x)>\ell(x)$ can not happen. Otherwise $\gamma_1, \cdots, \gamma_m, \alpha$ is a longer sequence satisfying the condition in \Cref{deflk}. This contradicts the maximality of $m$.
	
\end{proof}



\begin{example}
	Let $\Phi=E_6$ and $\lambda=\rho-2\alpha_1=(-1,2,3,4,5,-3,-3,3)\in\frh^*$. Here $\Delta=\{\alpha_i\}_{1\leq i\leq 6}$ with $\alpha_1=\frac{1}{2}(\ep_1+\ep_8-(\ep_2+\cdots+\ep_7))$, $\alpha_2=\ep_1+\ep_2$ and $\alpha_i=\ep_{i-1}-\ep_{i-2}$ for $3\leq i\leq 6$, where $\{\ep_i\}_{1\leq i\leq 8}$ is the  orthonormal basis of $\mathbb{R}^8$ (see \S\ref{subsec:RS}). Then $\langle\lambda, \alpha_1^\vee\rangle=-3$ and $\langle\lambda, \alpha_i^\vee\rangle>0$ for $i\neq 1$. Moreover, $s_{\alpha_1}\lambda=\rho+\alpha_1$ is dominant and semiregular ($\langle\rho+\alpha_1, \alpha_3^\vee\rangle=0$). By Theorem \ref{thmin1}, $L(\lambda)$ is integral minimal with $\gkd L(\lambda)=\langle\rho, \beta_s^{\vee}\rangle=11$. In view of Table \ref{tab1} and Theorem \ref{thmav}, it is minimal with associated variety $\overline{Be_{\alpha_1}}$. This is compatible with \cite[Ex. 7.2]{BXX}.
\end{example}

The proof of \Cref{thmin1} ``(2) implies (3)'' also leads to the following.

\begin{Cor}\label{corinw}
	Let $\lambda\in\frh^*$ be integral and $\nu$ be the dominant vector in $W\lambda$.  \begin{enumerate}
		\item If $\lambda$ is regular, then $L(\lambda)$ is minimal if and only if $\lambda=w\nu$ with $w\in\caC$.
		\item  If $\lambda$ is semiregular, then $L(\lambda)$ is minimal if and only if $\lambda=w\nu$ with $w$ or $ws_\alpha\in\caC$, where $\alpha$ is the unique simple root such that $\langle\nu, \alpha^\vee\rangle=0$. 
	\end{enumerate}
\end{Cor}

\begin{Rem}
	In \cite{Jo98}, Joseph classified all the minimal highest weight modules $L(\lambda)$ such that each $\Ann(L(\lambda))$ is completely prime. Most of these modules are presented in loc. cit. Table 3.
 One can verify that those cases are covered by the above corollary, while the nonintegral cases are contained in Theorem \ref{thmni} in the next section. Combined with Joseph's result, one can see that in most cases, $\Ann(L(\lambda))$ is not completely prime when $L(\lambda)$ is minimal.
\end{Rem}

Combining \Cref{thmav} and \Cref{thmin1} 
we have the following result. 
\begin{Cor}\label{corin}
 Suppose  $L(\lambda)$ is integral and minimal.
 Then 
 there is a unique simple root $\alpha$ such that $\inn{\lambda}{\alpha}\leq 0$ and 
 \[
 V(L(\lambda))=\overline{Be_{\alpha}}. 
 \]
\end{Cor}
Note that the above case only occurs when $\fgg$ has a simply laced root system. Moroever, $\alpha = \gamma_1$ in \Cref{deflk}. 

\trivial[h]{
There is some remarks: 
Obviously,  $\alpha = \gamma_1$ in \Cref{deflk}.
By \Cref{thmin1}, 
$w_0 w_{\lambda}^{-1} = s_{\gamma_m}\cdots, s_{\gamma_1} \in \cC$.
Note that $\cC$ decomposes into left/right cells according to the descent set \cite[Propositon~3.8]{L83}: For example, let $\cC_s := \set{w\in \cC | l(ws)<l(w)}$ and then $\cC_s$ are left cells and $\cC = \bigsqcup_{s\in S} \cC_s$.

Similarly, define right cell 
\[
{}_{s}\cC := \set{w\in \cC | l(sw)<l(w)}
\]

Now we see that $w_0w_\lambda^{-1} \in \cC_{s_{\gamma_1}}$. Note that taking inverse sends left/right cell to right/left cell. 
so $w_\lambda w_0 \in {}_{s_{\gamma_1}}\cC$. 
Now multiply $w_0$ on left or right send left/right cell to left/right cell. 
We conclude that $w_\lambda \sim_R s_{\gamma_i} w_0$

}

%
%
%
%
%
%
%

The following proposition provides criteria of $\nu$-minimal highest weight modules.

\begin{Pro}\label{proin1}
Suppose $\nu\in\mathfrak{D}$ is integral and non-regular.  For $\lambda\in W\nu$, the following conditions are equivalent:
	\begin{enumerate}
		\item  $\gkd L(\lambda)=|\Phi^+|-\aff(w_0w_{J_{w_0\nu}})$.
		\item  $w_\lambda\sim_{L}w_0w_{J_{w_0\nu}}$ in $W$.
		\item  $L(\lambda)$ is $\nu$-minimal.
	\end{enumerate}
\end{Pro}
Here $J_{w_0 \nu}$ is defined by \eqref{eq:Jmu}.

\begin{proof}
	Since $\nu$ is non-regular, the set $J:= J_{w_0\nu} \neq \emptyset$ and thus $w_{J}\neq 1$.
	
	Recall that $w_\lambda$ is a shortest representative in $W/W_{J}$. There exists $x\in W$ such that $xw_\lambda w_J=w_0$ and $\ell(x)+\ell(w_\lambda)+\ell(w_J)=\ell(w_0)$. In view of Lemma \ref{hklem1}, one has $w_0w_J\leq_L w_\lambda$ since $xw_{\lambda}=w_0w_J$ and $\ell(x)+\ell(w_\lambda)=\ell(w_0)-\ell(w_J)=\ell(w_0w_J)$. It follows from \Cref{alem1}(2) that $\aff(w_0w_J)\geq \aff(w_\lambda)$. Combined with \Cref{pr:main1} and \Cref{alem2}, one has 
	\[
	\gkd L(\lambda)=|\Phi^+|-\aff(w_\lambda)\geq|\Phi^+|-\aff(w_0w_J)>0. 
	\]
	By \cite[\S14.2 P9]{Lus03}, the equality holds if and only if $w_\lambda\sim_{L}w_0w_J$, 
 i.e., 
 $L(\lambda)$ is $\nu$-minimal.
	
\end{proof}

\begin{Rem}
    Fixing $\nu$ in \Cref{proin1}, there could be more than one $\lambda\in W \nu$  such that $L(\lambda)$ is $\nu$-minimal.
When $\nu$ is the infinitesimal character attached to an even nilpotent orbit ${}^L \cO$ in the Langlands dual of $\fgg$,
a classical result of Barbasch-Vogan \cite[Proposition~5.28 and Corollary~5.29]{BV85} shows that the number of $\lambda$ such that $L(\lambda)$ is $\nu$-minimal equals the number of conjugacy classes of the Lusztig's canonical quotient  attached to the orbit ${}^L\cO$ (one can use \cite[Theorem~5.9]{BG} to translate Barbasch-Vogan's result to category $\cO$).
\end{Rem}

\section{Nonintegral minimal highest weight modules}

%
%

We call a highest weight module $M$ is \textit{nonintegral minimal} if $\gkd M$ is minimal among all nonintegral highest weight modules. 

Recall that each $\lambda\in\frh^*$ gives an integral subsystem $\Phi_{[\lambda]}$ with simple system $\Delta_{[\lambda]}$ and Weyl group $W_{[\lambda]}$.
In this section, we will prove the following result about nonintegral minimal modules.

\begin{Thm}\label{thmni}
  Let $\Phi$ be irreducible and $\lambda\in\frh^*$. Then $L(\lambda)$ is
  nonintegral minimal if and only if $\lambda$ is dominant regular and
  $\Phi_{[\lambda]}$ is given in Table $\ref{tab2}$.
  \begin{table}[h]
    \centering \renewcommand{\arraystretch}{1.5} \setlength\tabcolsep{5pt}
    \begin{tabular}{c|ccccccccc}
      \hline $ \Phi$ & $A_{n}$ & $B_n$ & $C_n$ & $ D_n $ & $ E_6 $ & $ E_7 $ & $ E_8 $ & $ F_4 $ & $ G_2 $ \\ \hline
      $\Phi_{[\lambda]}$    & $A_{n-1}$ & $B_1\times B_{n-1}$ & $D_n$ & $ D_{n-1} $ & $ D_5 $ & $ E_6 $ & $ A_1\times E_7 $ & $ C_4 $ & $ A_2 $      \\ \hline
    \end{tabular}
    \bigskip
    \caption{Maximal proper integral subsystems}
    \label{tab2}
  \end{table}

  In this case, $\gkd L(\lambda)=|\Phi^+|-|\Phi_{[\lambda]}^+|$. Moreover, a
  nonintegral minimal module is minimal when $\Phi$ is of type $A$, $B$, $C$,
  $F$ or $G$.
\end{Thm}



\begin{example}
	Let $\Phi=F_4$ with $\Delta=\{\ep_2-\ep_3, \ep_3-\ep_4, \ep_4, \frac{1}{2}(\ep_1-\ep_2-\ep_3-\ep_4)\}$ (see \S\ref{subsec:RS}). Choose $\lambda=(4, 1, \frac{3}{2}, \frac{1}{2})$. It is easy to verify that $\Phi_{[\lambda]}$ is a subsystem with simple roots $\{\ep_3-\ep_4, \ep_4, \frac{1}{2}(\ep_1-\ep_2-\ep_3-\ep_4), \ep_2\}$. Therefore $\Phi_{[\lambda]}\simeq C_4$ and $\lambda$ is dominant. By Theorem \ref{thmni}, $L(\lambda)$ is nonintegral minimal with $\gkd L(\lambda)=|\Phi^+|-|\Phi_{[\lambda]}|=24-16=8$. It is minimal with associated variety $\overline{Be_{\ep_2-\ep_3}}$ in view of Table \ref{tab1} and Theorem \ref{thmav}. 
\end{example}


For the proof of Theorem \ref{thmni}, we need to recall some definitions and properties about closed subsystems and maximal subsystems. Most of which can be found in \cite[\S12]{Kan01}.

\begin{definition}\label{defcl}
	We say a subsystem $\Phi'\subset\Phi$ is \textit{closed} if $\alpha+\beta\in\Phi'$ whenever $\alpha, \beta\in\Phi'$ and $\alpha+\beta\in\Phi$.
\end{definition}

\begin{example}
	Let $\Phi=B_2=\{\pm (\ep_1+\ep_2), \pm (\ep_1-\ep_2), \pm \ep_1, \pm \ep_2\}$, then $\Phi'=\{\pm \ep_1, \pm \ep_2\}\simeq A_1\times A_1$ form a subsystem, which is not closed. 
\end{example}

The following result is an immediate consequence of Definition \ref{defcl}.

\begin{Lem}\label{lemni1}
	Let $\lambda\in\frh^*$. Then $\Phi_{[\lambda]}^\vee$ is a closed subsystem of $\Phi^\vee$.
\end{Lem}

Let $\alpha_0$ be the highest root of $\Phi$ and write
\[
\alpha_0=\sum_{i=1}^nh_i\alpha_i
\]
where $h_i$ are non-negative integers. 

\begin{Thm}[Borel-de Siebenthal]\label{thmbd}
	Let $\Phi$ be irreducible. Then the maximal proper closed subsystem of $\Phi$ $($up to the action of $W$$)$ are those with simple systems
	\begin{itemize}
		\item [(1)] $\{\alpha_1, \alpha_2, \cdots, \hat\alpha_i, \cdots, \alpha_n\}$ with $h_i=1$.
		\item [(2)] $\{-\alpha_0, \alpha_1, \cdots, \hat\alpha_i, \cdots, \alpha_n\}$ with $h_i$ 
  being a prime number.
	\end{itemize}
\end{Thm}

\begin{Pro}\label{proms}
  Let $\Phi$ be irreducible and $\lambda\in\frh^*$ be nonintegral.
  Then $|\Phi_{[\lambda]}|$ is maximal
  if and only $\Phi_{[\lambda]}$ has the type given by Table $\ref{tab2}$.

\end{Pro}

\begin{proof}
  We only give the proof when $\Phi=B_n$ since the other cases can be proved by a
  similar argument.
  %

  In this case $\Phi^\vee\simeq C_n$ (see \S\ref{subsec:RS}) with
  simple roots $\{\ep_1-\ep_2, \cdots, \ep_{n-1}-\ep_n, 2\ep_n\}$.
  The highest root is $\alpha_0=2\ep_1$.
  By \Cref{lemni1}, $\Phi_{[\lambda]}^\vee$ is a proper closed subsystem of $\Phi^{\vee}$. Since $|\Phi_{[\lambda]}| = |\Phi_{[\lambda]}^\vee|$ is maximal,   $\Phi_{[\lambda]}^\vee$ must be a maximal proper closed subsystem of $\Phi^\vee \cong  C_i\times C_{n-i}$ in view of \Cref{thmbd}. 
  Hence $|\Phi_{[\lambda]}|$ is maximal if and only if
  $\Phi_{[\lambda]}$ has type $ C_1^\vee\times C_{n-1}^\vee =  B_1\times B_{n-1}$.
  \end{proof}


\noindent\textbf{Proof of Theorem \ref{thmni}.} By Proposition \ref{pr:main1}, one has $\gkd L(\lambda)=|\Phi^+|-\aff_{[\lambda]}(w_\lambda)$. Lemma \ref{alem2} implies that $\aff_{[\lambda]}(w_\lambda)$ achieves its maximal value $|\Phi_{[\lambda]}^+|$ if and only if $w_\lambda$ is longest in $W_{[\lambda]}$, which means $\lambda$ is dominant regular. At last, we apply Proposition \ref{proms} and get Table \ref{tab2}. Combined with Table \ref{tab1}, we can easily find that nonintegral minimal module is minimal if $\Phi$ is of type $A$, $B$, $C$, $F$ or $G$. \qed

From the above arguments and the proof in Proposition \ref{proin1}, we have the following result.
\begin{Cor}
  Let $w_{[\lambda]}$ be the longest element in $W_{[\lambda]}$. Let
  $\nu\in\mathfrak{D}$ and $\lambda\in W\nu$. Let $\delta$ be the unique
  dominant weight in $W_{[\lambda]}\lambda$. If $\nu$ is singular or
  nonintegral, then the following conditions are equivalent:
  \begin{enumerate}
    \item  $L(\lambda)$ is $\nu$-minimal.
    \item
          $\gkd L(\lambda)=\gkd L(\delta)=|\Phi^+|-\aff_{[\lambda]}(w_{[\lambda]}w_{J_\lambda})$.
    \item  $w_\lambda\sim_{L}w_{[\lambda]}w_{J_\lambda}$ in $W_{[\lambda]}$.
  \end{enumerate}
\end{Cor}

\trivial[h]{
\section{Minimal highest weight modules of classical types}\label{classical}

%
%
In order to verify whether a highest weight module $L(\lambda)$ is minimal, we apply
Theorem \ref{thmin1} for types $A$, $D$, $E$ when $\lambda$ is integral and Theorem \ref{thmni} for types $A$, $B$, $C$, $F$, $G$ when $\lambda$ is nonintegral. These criteria are quite practical. In this section, we intend to investigate more explicit characterization of minimal highest weight modules based on these criteria in the case of classical Lie algebras.

We continue to follow the notation in \cite{Hum78}. Note that for $\Phi=A_{n-1}, B_n, C_n$ and $D_n$, a weight $ \lambda\in \hs $ will be identified with a sequence 
$ \lambda=(\lambda_1,\cdots,\lambda_n)\in\mathbb{C}^n $
such that $ \lambda=\sum_{i=1}^n\lambda_i\varepsilon_i $.

\begin{definition}
	Following Mathieu \cite{Ma}, a sequence $ x=(x_1,x_2,\cdots, x_m) \in \mathbb{C}^m$ is said to be \emph{ordered} if all differences $x_i-x_{i+1}$ are positive integers.
\end{definition}

\subsection{Type $ A_{n} $}

\begin{Thm}\label{type-A}
	Let $\fgg=\mathfrak{sl}(n+1, \mathbb{C})$. The simple module $ L(\lambda)$ is  minimal  if and only if the  length of the longest ordered subsequence of $\lambda$ is $ n $. In this case, let $t$ be the smallest integer such that $(\lambda_1, \cdots, \hat{\lambda}_t, \cdots, \lambda_{n+1})$ is ordered. Then
	\begin{itemize}
		\item [(1)] If $\lambda$ is not integral, then $V(L(\lambda))=\overline{Be_{\alpha_{t-1}}}\cup \overline{Be_{\alpha_{t}}}$, where $e_{\alpha_0}=e_{\alpha_{n+1}}=0$.
		\item [(2)] If $\lambda$ is integral, then $V(L(\lambda))=\overline{Be_{\alpha_{t-1}}}$ when $\lambda_{t-1}\leq \lambda_t$ and $V(L(\lambda))=\overline{Be_{\alpha_{t}}}$ when $\lambda_{t}\leq \lambda_{t+1}$.
	\end{itemize}

\end{Thm}
\begin{proof}
	It suffices to consider the first assertion since results about associated varieties are evident consequences of Theorem \ref{thmav}.
	
	If $\lambda$ is nonintegral, the first statement follows from Theorem \ref{thmni}. Indeed, the module $L(\lambda)$ is minimal if and only if $\lambda$ is regular dominant and $\Phi_{[\lambda]}\simeq A_{n-1}$, which is equivalent to the required ordered subsequence. 
	
	Now suppose that $\lambda$ is integral. Set $w=w_0w_\lambda^{-1}$. It is the longest element such that $w\lambda$ is dominant. If $L(\lambda)$ is minimal, by Corollary \ref{corinw}, we can assume that $w\in\mathcal{C}$, which means $w\neq1$ has the unique reduced expression. So there exist some $i, m\geq1$ such that (recall that $\alpha_j=\ep_j-\ep_{j+1}$ for $1\leq j\leq n$)
	\[
	w=s_{\alpha_{i+m-1}}\cdots s_{\alpha_{i+1}}s_{\alpha_i}\ \mbox{or}\ s_{\alpha_{i-m+1}} \cdots s_{\alpha_{i-1}}s_{\alpha_{i}}.
	\]
	Since the arguments are similar, we only consider the former case, for which
	\[
	w\lambda=(\lambda_1, \cdots, \lambda_{i-1}, \lambda_{i+1}, \cdots, \lambda_{i+m}, \lambda_i, \lambda_{i+m+1}, \cdots, \lambda_{n+1})
	\]
	is dominant and thus $\lambda_{i-1}\geq\lambda_{i+1}$. It follows from Theorem \ref{thmin1}(3) ($k=1$ and $\gamma_1=\alpha_i$) that $\langle\lambda, \alpha_j^\vee\rangle=\lambda_j-\lambda_{j+1}\in\mathbb{Z}_{>0}$ for $j\neq i$. If $\lambda_{i-1}>\lambda_{i+1}$, then $\lambda$ is ordered by removing $\lambda_i$. If $\lambda_{i-1}=\lambda_{i+1}$, then $\langle s_{\alpha_i}\lambda, \alpha_{i-1}^\vee\rangle=0$. Theorem \ref{thmin1}(3) yields $\gamma_m=\alpha_{i+m-1}=\alpha_{i-1}$, which is impossible.
	
	Conversely, assume that $(\lambda_1, \cdots, \hat{\lambda}_t, \cdots, \lambda_{n+1})$ is ordered and $\lambda$ is not, then $(\lambda_{t-1}-\lambda_t)(\lambda_t-\lambda_{t+1})\leq0$. If $\lambda_t\geq\lambda_{t-1}$, there exists $m\geq1$ such that $\lambda_{t-m-1}>\lambda_t\geq\lambda_{t-m}$ (set $\lambda_0=+\infty$). Evidently
	\[
	(\lambda_1, \cdots, \lambda_{t-m-1}, \lambda_t, \lambda_{t-m}, \cdots, \lambda_{t-1}, \lambda_{t+1}, \cdots, \lambda_{n+1})=s_{\alpha_{t-m}}\cdots s_{\alpha_{t-1}}\lambda
	\]
	is dominant and $w=s_{\alpha_{t-m}}\cdots s_{\alpha_{t-1}}\in\mathcal{C}$. Therefore $L(\lambda)$ is minimal by Corollary \ref{corinw}. If $\lambda_t\leq\lambda_{t+1}$, the argument is similar.
	
\end{proof}

%



\subsection{Type $B_n$}

\begin{Thm}\label{type-B}
	Let $\fgg=\mathfrak{so}(2n+1, \mathbb{C})$. The simple module $L(\lambda)$ is minimal if and only if the following conditions are satisfied:
	\begin{itemize}
		\item [(1)] There is $1\leq t\leq n$ such that $(\lambda_1, \cdots, \lambda_n)$ is ordered after $\lambda_t$ is removed.
		\item [(2)] $2\lambda_i\in\mathbb{Z}_{>0}$ for $1\leq i\leq n$, while $2\lambda_t$ and $2\lambda_i$ $(i\neq t)$ have different parity.
	\end{itemize}
In this case, $V(L(\lambda))=\overline{Be_{\alpha_{t-1}}}\cup \overline{Be_{\alpha_{t}}}$ for $t<n$ and  $V(L(\lambda))=\overline{Be_{\alpha_{n-1}}}$ for $t=n$.
\end{Thm}
\begin{proof}
By Theorem \ref{thmin1} and \ref{thmni}, $L(\lambda)$ is minimal if and only if $\lambda$ is nonintegral dominant regular and $\Phi_{[\lambda]}\simeq B_1\times B_{n-1}$. This yields (1) and (2). The last assertion follows from Theorem \ref{thmav}.
\end{proof}

%

\subsection{Type $C_n$}

%

\begin{Thm}\label{cthm}
	Let $\fgg=\mathfrak{sp}(n, \mathbb{C})$. The simple module $L(\lambda)$ is minimal if and only if the following conditions are satisfied:
	\begin{itemize}
		\item [(1)] $\lambda_i\in\frac{1}{2}+\mathbb{Z}$ for all $1\leq i\leq n$.
		\item [(2)] $(\lambda_1, \cdots, \lambda_{n-1}, |\lambda_n|)$ is ordered.
	\end{itemize}
	In this case, $V(L(\lambda))=\overline{Be_{\alpha_{n}}}$.
\end{Thm}

\begin{proof}
	This follows easily from Theorem \ref{thmni} and \ref{thmav}.
\end{proof}

%
%

\subsection{Type $D_n$}

\begin{Thm}\label{dthm}
	Let $\Phi=D_{n}$ and $\lambda\in\frh^*$. Then $L(\lambda)$ is minimal if and only if the following conditions are satisfied:
	\begin{itemize}
		\item [(1)] $\lambda$ is integral and there is the unique $1\leq k\leq n$ such that $\langle\lambda, \alpha_k^\vee\rangle\leq0$.
		\item [(2)] If $1\leq k\leq n-2$, then $(\lambda_1, \cdots, \hat\lambda_l, \cdots \lambda_{n-1}, |\lambda_n|)$ is ordered and $\lambda_l>-|\lambda_n|$ for $l=k$ or $k+1$.
		\item [(3)] If $k=n-1, n$, then $\lambda_{n-2}>|\lambda_{n-1}|$ or $$\lambda_{n-3}>(-1)^{n-k+1}\lambda_n>-\lambda_{n-1}>-\lambda_{n-2}\geq\lambda_{n-1};$$
	\end{itemize} 
	In this case, $V(L(\lambda))=\overline{Be_{\alpha_{k}}}$.
\end{Thm}
\begin{proof}
	Set $w=w_0w_\lambda^{-1}$. Then $w$ is the longest element such that $w\lambda$ is dominant.
	
	By Theorem \ref{thmin1} and \ref{thmni}, $L(\lambda)$ is minimal if and only if $\lambda$ is integral and $w\in\mathcal{C}$. We list all the elements in $\mathcal{C}$.
	\smallskip
	
	(a) $w=s_{\alpha_{k+m-1}}\cdots s_{\alpha_{k}}$ with $k\leq n-2$ and $m<n-k$.
	
	(b) $w=s_{\alpha_{n-1}}s_{\alpha_{n-2}}\cdots s_{\alpha_{k}}$ with $k\leq n-2$ and $m=n-k$.
	
	(c) $w=s_{\alpha_{n}}s_{\alpha_{n-2}}\cdots s_{\alpha_{k}}$ with $k\leq n-2$ and $m=n-k$.
	
	(d) $w=s_{\alpha_{k-m+1}}\cdots s_{\alpha_{k}}$ with $k\leq n-2$ and $m\leq k$.
	
	(e) $w=s_{\alpha_{n-m}}\cdots s_{\alpha_{n-2}}s_{\alpha_{n-1}}$ with $k=n-1$ and $m<n$.
	
	(f) $w=s_{\alpha_{n-m}}\cdots s_{\alpha_{n-2}}s_{\alpha_{n}}$ with $k=n$ and $m<n$.
	
	(g) $w=s_{\alpha_{n}}s_{\alpha_{n-2}}s_{\alpha_{n-1}}$ with $k=n-1$.
	
	(h) $w=s_{\alpha_{n-1}}s_{\alpha_{n-2}}s_{\alpha_{n}}$ with $k=n$.
	\smallskip
	
	This makes a case-by-case proof possible. First assume that (1) and (2) are true with $l=k$ in (2). Then (1) and (2) yield $\lambda_{k+1}\geq\lambda_k>-|\lambda_n|$ and 
	\begin{equation*}\label{ineq0}
	\lambda_1>\cdots>\lambda_{k-1}>\lambda_{k+1}>\cdots>\lambda_{n-1}>|\lambda_n|.
	\end{equation*}
	Let $m$ be the largest integer such that $\lambda_{k+m}\geq\lambda_k$ (thus $\lambda_k>\lambda_{k+m+1}$ when $k+m+1\leq n$). If $\lambda_k>|\lambda_{n}|$ (thus $k+m<n$), it can be easily verified that $\alpha_{k+m-1}, \cdots, \alpha_{k}$ is the linked sequence of $\lambda$. Moreover
	\begin{equation}\label{ineq1}
		s_{\alpha_{k+m-1}}\cdots s_{\alpha_{k}}\lambda=(\lambda_1, \cdots, \lambda_{k-1}, \lambda_{k+1}, \cdots, \lambda_{k+m}, \lambda_k, \lambda_{k+m+1}, \cdots, \lambda_n) 
     \end{equation}
is dominant and regular or semiregular (when $\lambda_{k+m}=\lambda_k$). By Theorem \ref{thmin1}, this gives the case (a). If $\lambda_k\leq|\lambda_n|$ and $\lambda_n>0$, then $\alpha_{n-1}, \cdots, \alpha_{k}$ is the linked sequence of $\lambda$. In addition, 
	\begin{equation}\label{ineq2}
		s_{\alpha_{n-1}}s_{\alpha_{n-2}}\cdots s_{\alpha_{k}}\lambda=(\lambda_1, \cdots, \lambda_{k-1}, \lambda_{k+1}, \cdots, \lambda_{n-1}, \lambda_n, \lambda_k) 
	\end{equation}
	is dominant and regular or semiregular (when $\lambda_{k}=\lambda_n$). This is the case (b) by Theorem \ref{thmin1}. If $\lambda_k\leq|\lambda_n|$ with $\lambda_n\leq 0$, then $\alpha_n, \alpha_{n-2}, \cdots, \alpha_{k}$ is the linked sequence of $\lambda$ and 
	\begin{equation}\label{ineq3}
	s_{\alpha_{n}}s_{\alpha_{n-2}}\cdots s_{\alpha_{k}}\lambda=(\lambda_1, \cdots, \lambda_{k-1}, \lambda_{k+1}, \cdots, \lambda_{n-1}, -\lambda_n, -\lambda_k).
	\end{equation}
	We arrived at the case (c) ($\lambda$ is semiregular when $\lambda_k=-\lambda_n$). If (1) and (2) are true with $l=k+1$, then we get (d) for some $m\leq k$. Now suppose we have (1) and (3). Then $\lambda_{n-2}>|\lambda_{n-1}|$ gives (e) and (f). If $k=n-1$ and  $\lambda_{n-3}>\lambda_n>-\lambda_{n-1}>-\lambda_{n-2}\geq\lambda_{n-1}$, one has (g). If $k=n$ and  $\lambda_{n-3}>-\lambda_n>-\lambda_{n-1}>-\lambda_{n-2}\geq\lambda_{n-1}$, we obtain (h).
	
	Conversely, assume that $L(\lambda)$ is minimal and (a) holds. Then $\alpha_{k+m-1}, \dots,$ $ \alpha_{k}$ is the linked sequence of $\lambda$ and $s_{\alpha_{k+m-1}}\cdots s_{\alpha_{k}}\lambda$ is dominant, in view of Theorem \ref{thmin1}. Moreover, $\lambda$ is regular or semiregular. If $\lambda$ is regular, it is easy to see (1) and (2) are true for $l=k$, keeping in mind that the weight in (\ref{ineq1}) is dominant. If $\lambda$ is semiregular, the definition of linked sequence forces $\langle w\lambda, \alpha_{k+m-1}^\vee\rangle=0$, that is, $\lambda_k=\lambda_{k+m}$. With (\ref{ineq1}) in hand, (1) and (2) also hold. The arguments for other cases which we omit are similar.

\end{proof}

\mjjd{How about the exceptional case?}{}
}

\subsection*{Acknowledgments}
We would like to thank Toshiyuki Tanisaki for very helpful
conversations about orbital varieties.
We  would like to thank Binyong Sun for inspiring us this problem. We also  would like to thank the  referee for very helpful suggestions and comments.

\printbibliography

\end{document}